# GOOD ROUGH PATH SEQUENCES AND APPLICATIONS TO ANTICIPATING STOCHASTIC CALCULUS

By Laure Coutin, Peter Friz and Nicolas Victoir

*Université Paul Sabatier, University of Cambridge and University Oxford*

We consider anticipative Stratonovich stochastic differential equations driven by some stochastic process lifted to a rough path. Neither adaptedness of initial point and vector fields nor commuting conditions between vector field is assumed. Under a simple condition on the stochastic process, we show that the unique solution of the above SDE understood in the rough path sense is actually a Stratonovich solution. We then show that this condition is satisfied by the Brownian motion. As application, we obtain rather flexible results such as support theorems, large deviation principles and Wong–Zakai approximations for SDEs driven by Brownian motion along anticipating vectorfields. In particular, this unifies many results on anticipative SDEs.

**1. Introduction.** We fix a filtered probability space $(\Omega, \mathcal{F}, \mathbb{P}, (\mathcal{F}_t)_{0 \le t \le 1})$ satisfying the usual conditions. Itô's theory tells us that there exists a unique solution to the Stratonovich stochastic differential equation

$$(1.1) \qquad \begin{cases} dY_t = V_0(Y_t)\,dt + \sum_{i=1}^{d} V_i(Y_t) \circ dB_t^i, \\ Y_0 = y_0, \end{cases}$$

where $B = (B_t^1, \ldots, B_t^d)_{0 \le t \le 1}$ is a standard $d$-dimensional $(\mathcal{F}_t)_{0 \le t \le 1}$-Brownian motion, $V_0$ a $C^1$ vector field on $\mathbb{R}^d$, $V_1, \ldots, V_d$ some $C^2$ vector fields on $\mathbb{R}^d$, and $y_0$ a $\mathcal{F}_0$-measurable random variable.

We remind the reader that, by definition (see [19]), a process $z$ is Stratonovich integrable with respect to a process $x$ if for all $t$ there exists a random variable denoted $\int_0^t z_u \circ dx_u$ (the Stratonovich integral of $z$ with respect to $x$ between time 0 and $t$) such that for all sequences $(D^n = (t_i^n)_i)_{n \ge 0}$ of subdivisions of $[0, t]$ such that $|D^n| \to_{n \to \infty} 0$, the following convergence holds in

probability

$$\sum_i \left( \frac{1}{t_{i+1}^n - t_i^n} \int_{t_i^n}^{t_{i+1}^n} z_u \, du \right) (x_{t_{i+1}^n} - x_{t_i^n}) \underset{n \to \infty}{\to} \int_0^t z_u \circ dx_u.$$

Another way to express this is by introducing $x^D$ the $D$-linear approximation of $x$, where $D = (t_i)$ is a subdivision of $[0,1]$:

$$(1.2) \qquad x^D(t) = x_{t_i} + \frac{t - t_i}{t_{i+1} - t_i}(x_{t_{i+1}} - x_{t_i}) \qquad \text{if } t_i \le t \le t_{i+1}.$$

Then, $z$ is Stratonovich integrable with respect to $x$ if and only if for all sequences of subdivisions $D^n$ which mesh size tends to 0 and for all $t \in [0,1]$, $\int_0^t z_u \, dx^{D^n}(u)$ converges in probability to $\int_0^t z_u \circ dx_u$.

Ocone and Pardoux [20] showed that there exists a unique solution to equation (1.1) even if the vector field $V_0$ and the initial condition $y_0$ were allowed to be $(\mathcal{F}_1)$-random variables. They did so relating the Skorokhod integral and the Stratonovich one, and using Malliavin calculus techniques. This solution has been studied in various directions: existence and study of the density of $Y$ [2, 22, 23], a Freidlin–Wentzell-type theorem [18], results on the support of the law of $\{Y_t, 0 \le t \le 1\}$ [3, 17], approximation of $Y_t$ by some Euler's type schemes [1].... The case where the vector fields $V_1, \ldots, V_d$ are allowed to be $(\mathcal{F}_1)$-measurable was dealt in [8, 9], under the strong condition that the $(V_i)_{1 \le i \le d}$ commute. This is an application of the Doss–Sussman theorem. The latter says that if $V = (V_1, \ldots, V_d)$ are $d$ vector fields smooth enough such that $[V_i, V_j] = 0$, for all $i, j \ge 1$, and if $V_0$ is another vector fields, then the map $\varphi_{(V_0,V)}^{y_0}$ which at a smooth path $x : [0,1] \to \mathbb{R}^d$ associates the path $y$ which is the solution of the differential equation

$$(1.3) \qquad \begin{cases} dy_t = V_0(y_t) \, dt + \sum_{i=1}^d V_i(y_t) \, dx_t^i, \\ y_0 = y_0, \end{cases}$$

is continuous when one equips the space of continuous functions with the uniform topology. One can then define $\varphi_{(V_0,V)}^{y_0}$ on the whole space of continuous function, in particular $\varphi_{(V_0,V)}^{y_0}(B)$ is then well defined, and is almost surely the solution of the Stratonovich differential equation (1.1). This remains true even if the vector fields and the initial condition are allowed to be random.

Rough path theory can be seen as a major extension of the Doss–Sussman result. One of the main thing to remember from this theory is that it is not $x$ which controls the differential equation (1.3), but the lift of $x$ to a path in a Lie group lying over $\mathbb{R}^d$. The choice of the Lie group depends on the roughness of $x$.

If $x$ is a $\mathbb{R}^d$-valued path of finite $p$-variation, $p \geq 1$, one needs to lift $x$ to a path $\mathbf{x}$ with values in $G^{[p]}(\mathbb{R}^d)$, the free nilpotent Lie group of step $[p]$ over $\mathbb{R}^d$. When $x$ is smooth, there exists a canonical lift of $x$ to a path denoted $S(x)$ with values in $G^{[p]}(\mathbb{R}^d)$ ($S(x)$ is obtain from $x$ by computing the "first $[p]$" iterated integrals of $x$). If $x$ is a smooth path, then there exists a solution $y$ to the differential equation (1.3). We denote by $I^{y_0}_{(V_0, V)}$ the map which at $S(x)$ associates $S(x \oplus y)$. Denoting $C(G^{[p]}(\mathbb{R}^d))$ the set of continuous paths from $[0, 1]$ into $G^{[p]}(\mathbb{R}^d)$, we see that $I^{y_0}_{(V_0, V)}$ is a map from a subset of $C(G^{[p]}(\mathbb{R}^d))$ to $C(G^{[p]}(\mathbb{R}^d \oplus \mathbb{R}^n))$. Lyons showed that this map is (locally uniformly) continuous when one equips $C(G^{[p]}(\mathbb{R}^d))$ and $C(G^{[p]}(\mathbb{R}^d \oplus \mathbb{R}^n))$ with a "$p$-variation distance." Hence one can define $I^{y_0}_{(V_0, V)}$ on the closure (in this $p$-variation topology) of the canonical lift of smooth paths. This latter set is the set of geometric $p$-rough paths.

In the case $x = B$, the Brownian motion being almost surely $1/p$-Hölder, $2 < p < 3$, one needs to lift $B$ to a process with values in $G^2(\mathbb{R}^d)$ to obtain the solution in the rough path sense of equation (1.1). This is equivalent to define its area process. The standard choice for the area process of the Brownian motion is the Lévy area [13, 15], although one could choose very different area processes [11]. Choosing this area, we lift $B$ to a geometric $p$-rough path $\mathbf{B}$, $I^{y_0}_{(V_0, V)}(\mathbf{B})$ is then the Stratonovich solution of the stochastic differential equation (1.1), together with its lift (i.e. here its area process).

Just as before, when the vector fields $V_i$ are almost surely "smooth enough" and $y_0$ is almost surely finite, the Itô map $I^{y_0}_{(V_0, V)}$ is still well defined and continuous almost surely, and there is no problem at all of definition of $I^{y_0}_{(V_0, V)}(\mathbf{B})$. Therefore, the theory of rough path provides a meaning and a unique solution to the stochastic differential equation (1.1), even when the vector fields and the initial condition depend on the whole Brownian path. Moreover, the continuity of the Itô map provides for free a Wong–Zakai theorem, and is very well adapted to obtaining large deviation principles and support theorems.

The only work not completely given for free by the theory of rough path is to prove that the solution $\mathbf{y}$ of equation (1.1) using the rough path approach is actually solution of the Stratonovich differential equation, that is, that for all $t$,

$$y_t = y_0 + \int_0^t V_0(y_u)\, du + \sum_{i=1}^d \int_0^t V_i(y_u) \circ dB_u^i.$$

When the vector fields and the initial condition are deterministic, this is usually proved using the standard Wong–Zakai theorem.

We provide in this general case here a solution typically in the spirit of rough path, by separating neatly probability theory and differential equation theory. We will show via a deterministic argument that to obtain our result we only need to check that, if $D^n$ is a sequence of subdivisions which steps tends to 0, $\int_0^\cdot B_u \otimes dB_u^{D^n}$ and $\int_0^\cdot B_u^{D^n} \otimes dB_u^{D^n}$ converges in an appropriate topology to $\int_0^\cdot B_u \otimes \circ dB_u$.

The paper is organized as follows: in the first section, we present quickly the theory of rough path (for a more complete presentation; see [14, 15] or [10]). The second section introduce the notion of good rough path sequence and its properties. We will then show that $B^n$ defines a good rough path sequence associated to **B**, and this will imply that the solution via rough path of equation (1.3) with signal **B** is indeed solution of the Stratonovich stochastic differential equation (1.1).

We conclude with a few applications: a Wong–Zakai-type approximation result, a large deviation principle, and some remarks on the support theorem.

**2. Rough paths.** By path we will always mean a continuous function from $[0,1]$ into a (Lie) group. If $x$ is such a path, $x_{s,t}$ is a notation for $x_s^{-1} \cdot x_t$.

2.1. *Algebraic preliminaries.* We present the theory of rough paths, in the finite dimensional case purely for simplicity. All arguments are valid in infinite dimension. We equip $\mathbb{R}^d$ with the Euclidean scalar product $\langle \cdot, \cdot \rangle$ and the Euclidean norm $|x| = \langle x, x \rangle^{1/2}$. We denote by $(G^n(\mathbb{R}^d), \otimes)$ the free nilpotent of step $n$ over $\mathbb{R}^d$, which is imbedded in the tensor algebra $T^n(\mathbb{R}^d) = \bigoplus_{k=0}^n (\mathbb{R}^d)^{\otimes k}$. We define a family of dilations on the group $G^n(\mathbb{R}^d)$ by the formula

$$\delta_\lambda(1, x_1, \ldots, x_n) = (1, \lambda x_1, \ldots, \lambda^n x_n),$$

where $x_i \in (\mathbb{R}^d)^{\otimes i}$, $(1, x_1, \ldots, x_n) \in G^n(\mathbb{R}^d)$ and $\lambda \in \mathbb{R}$. Inverse, exponential and logarithm functions on $T^n(\mathbb{R}^d)$ are defined by mean of their power series [14, 21]. We define on $(\mathbb{R}^d)^{\otimes k}$ the Hilbert tensor scalar product and its norm

$$\left\langle \sum_{i=1}^d x_i^1 \otimes \cdots \otimes x_i^k, \sum_{j=1}^d y_j^1 \otimes \cdots \otimes y_j^k \right\rangle = \sum_{i,j} \langle x_i^1, y_j^1 \rangle \cdots \langle x_i^k, y_j^k \rangle,$$

$$|x| = \langle x, x \rangle^{1/2}.$$

This yields a family of compatible tensor norms on $\mathbb{R}^d$ and its tensor product spaces. Since all finite dimensional norms are equivalent the Hilbert structure of $\mathbb{R}^d$ was only used for convenience. In fact, one can replace $\mathbb{R}^d$ by a Banach space and deal with suited tensor norms but this can be rather subtle; see [15] and the references therein.

For $(1, x_1, \ldots, x_n) \in G^n(\mathbb{R}^d)$, with $x_i \in (\mathbb{R}^d)^{\otimes i}$, we define

$$\|(1, x_1, \ldots, x_n)\| = \max_{i=1,\ldots,n} \{(i! |x_i|)^{1/i}\}.$$

$\|\cdot\|$ is then a symmetric subadditive homogeneous norm on $G^n(\mathbb{R}^d)$ ($\|g\| = 0$ iff $g = 1$, and $\|\delta_\lambda g\| = |\lambda| \cdot \|g\|$ for all $(\lambda, g) \in \mathbb{R} \times G^n(\mathbb{R}^d)$, for all $g, h \in G^n(\mathbb{R}^d)$, $\|g \otimes h\| \le \|g\| + \|h\|$, and $\|g^{-1}\| = \|g\|$). Of course, we could have chosen any other continuous homogeneous norm, as they are all equivalent (see [15]). This homogeneous norm allows us to define on $G^n(\mathbb{R}^d)$ a left-invariant distance with the formula

$$d(g, h) = \|g^{-1} \otimes h\|.$$

From this distance we can define some distances on the space of continuous paths from $[0, 1]$ into $G^{[p]}(\mathbb{R}^d)$ ($p$ is a fixed real greater than or equal to 1):

(1) The $p$-variation distance

$$d_{p\text{-var}}(x, y) = d(x_0, y_0) + \sup \left( \sum_i d(x_{t_i, t_{i+1}}, y_{t_i, t_{i+1}})^p \right)^{1/p}$$

where the supremum is over all subdivision $(t_i)_i$ of $[0, 1]$. We also define $\|x\|_{p\text{-var}} = d_{p\text{-var}}(1, x)$, and $\|x\|_{p\text{-var}} < \infty$ means that $x$ has finite $p$-variation.

(2) Modulus distances: We say that $\omega \colon \{(s, t), 0 \le s \le t \le 1\} \to \mathbb{R}^+$ is a control if

$$\begin{cases} \omega \text{ is continuous.} \\ \omega \text{ is super-additive, that is } \forall s < t < u, \omega(s, t) + \omega(t, u) \le \omega(t, u). \\ \omega \text{ is zero on the diagonal, that is } \omega(t, t) = 0 \text{ for all } t \in [0, 1]. \end{cases}$$

When $\omega$ is non-zero off the diagonal we introduce the distance $d_{\omega, p}$ by

$$d_{\omega, p}(x, y) = d(x_0, y_0) + \sup_{0 \le s < t \le 1} \frac{d(x_{s,t}, y_{s,t})}{\omega(s, t)^{1/p}}.$$

The simplest example of such a control is given by $\omega(s, t) = t - s$. In this case, $d_{\omega, p}$ defines then a notion of $1/p$-Hölder distance on $G^n(\mathbb{R}^d)$-valued paths. In general, if $\|x\|_{\omega, p} = d_{\omega, p}(1, x) < \infty$, we say that $x$ has finite $p$-variation controlled by $\omega$.

2.2. *The $G^n(\mathbb{R}^d)$-valued paths.* Let $x \colon [0, 1] \to \mathbb{R}^d$ be a path of bounded variation and define $S(x)$ to be the solution of the ordinary differential equation

$$dS(x)_t = S(x)_t \otimes dx_t,$$

$$S(x)_0 = \exp(x_0),$$

where $\otimes$ is the multiplication of $T^n(\mathbb{R}^d)$ and exp is the exponential on $T^n(\mathbb{R}^d)$. $S(x)$ actually takes its values in $G^n(\mathbb{R}^d)$. Moreover, if the 1-variation of $x$ is controlled by $\omega$, then so is the 1-variation of $S(x)$ ([14], Theorem 1).

EXAMPLE 1.   If $n = 2$, $S(x)_t = \exp(x(t) + \frac{1}{2}\int_0^t x_u \otimes dx_u - \frac{1}{2}\int_0^t dx_u \otimes x_u)$. The term $\frac{1}{2}\int_0^t x_u \otimes dx_u - \frac{1}{2}\int_0^t dx_u \otimes x_u$ is the area between the line joining $x_0$ and $x_t$ and the path $(x_u)_{0 \le u \le t}$.

DEFINITION 1.   A path $\mathbf{x}\colon [0,1] \to G^{[p]}(\mathbb{R}^d)$ is a geometric $p$-rough path if there exists a sequence of paths of bounded variation $(x_n)_{n \in \mathbb{N}}$ such that

$$\lim_{n \to \infty} d_{p\text{-var}}(S(x_n), \mathbf{x}) = 0.$$

NOTATION 1.   Let $C^{0,p\text{-var}}(G^{[p]}(\mathbb{R}^d))$ denote the space of geometric $p$-rough paths, and $C^{0,\omega,p}(G^{[p]}(\mathbb{R}^d))$ the set of paths $\mathbf{x}\colon [0,1] \to G^{[p]}(\mathbb{R}^d)$ for which there exists a sequence of paths $(x_n)_{n \in \mathbb{N}}$ of bounded variation such that

$$\lim_{n \to \infty} d_{\omega,p}(S(x_n), \mathbf{x}) = 0.$$

In particular, any $G^{[p]}(\mathbb{R}^d)$-valued path with finite $q$-variation, $q < p$, is a geometric $p$-rough path, and a geometric $p$-rough path has finite $p$-variation. We refer to [7] for a "Polish" study of the space of geometric rough paths.

If $\mathbf{x}$ is a geometric $p$-rough path, we denote by $\mathbf{x}_{s,t}^i$ the projection of $\mathbf{x}$ onto $(\mathbb{R}^d)^{\otimes i}$. We also say that $\mathbf{x}$ lies above $\mathbf{x}^1$.

2.3. *The Itô map.*   We fix $p \ge 1$, and a control $\omega$.

When $A$ is vector field on $\mathbb{R}^d$, we denote by $d^k A$ its $k^{th}$ derivative (with the convention $d^0 A = A$), and its $\gamma$-Lipschitz norm by

$$\|A\|_{\mathrm{Lip}(\gamma)} = \max\left\{\max_{k=0,\dots,[\gamma]} \|d^k A\|_\infty, \|d^{[\gamma]} A\|_{\gamma - [\gamma]}\right\},$$

where $\|\cdot\|_\infty$ is the sup norm and $\|\cdot\|_\beta$ the $\beta$-Hölder norm, $0 \le \beta < 1$. If $\|A\|_{\mathrm{Lip}(\gamma)} < \infty$, we say that $A$ is a $\mathrm{Lip}(\gamma)$-vector fields on $\mathbb{R}^d$.

We consider $V = (V_1, \dots, V_d)$, where the $V_i$ are $\mathrm{Lip}(p + \varepsilon)$-vector fields on $\mathbb{R}^n$, $\varepsilon > 0$. $V$ can be identified with a linear map from $\mathbb{R}^d$ into $\mathrm{Lip}(p+\varepsilon)$-vector fields on $\mathbb{R}^n$,

$$V(y)(dx^1, \dots, dx^d) = \sum_{i=1}^d V_i(y)\,dx^i.$$

For a $\mathbb{R}^d$-valued path $x$ of bounded variation, we define $y$ to be the solution of the ordinary differential equation

$$(2.1) \qquad\qquad \begin{cases} dy_t = V(y_t)\,dx_t, \\ y_0 = y_0. \end{cases}$$

Lifting $x$ and $(x \oplus y)$ to $S(x)$ and $S(x \oplus y)$ (their canonical lift to paths with values in the free nilpotent group of step $[p]$), we consider the map which at $S(x)$ associates $S(x \oplus y)$. We denote it $I_{y_0,V}$. We refer to [14, 15] for the following theorem (which we adapt here to our notation, but its equivalence with the one stated in [14, 15] is a trivial consequence of Proposition 4 in [7]).

THEOREM 1 [Universal limit theorem (Lyons)]. *The map $I_{y_0,V}$ is continuous from $(C^{0,\omega,p}(G^{[p]}(\mathbb{R}^d)), d_{\omega,p})$ into $(C^{0,\omega,p}(G^{[p]}(\mathbb{R}^d \oplus \mathbb{R}^n)), d_{\omega,p})$.*

Let $x_n$ be a sequence of path of bounded variation such that $S(x_n)$ converges in the $d_{\omega,p}$-topology to a geometric $p$-rough path $\mathbf{x}$, and define $y_n$ to be the solution of the differential equation (2.1), where $x$ is replaced by $x_n$. Then, the universal limit theorem says that $S(x_n \oplus y_n)$ converges in the $d_{\omega,p}$-topology to a geometric $p$-rough path $\mathbf{z}$. We say that $\mathbf{y}$, the projection of $\mathbf{z}$ onto $G^{[p]}(\mathbb{R}^n)$ is the solution of the rough differential equation

$$d\mathbf{y} = V(\mathbf{y}) \, d\mathbf{x}$$

with initial condition $y_0$. It is interesting to observe that Lyons' estimates actually give that for all $R > 0$ and sequence $\mathbf{x}_n \in C^{0,\omega,p}(G^{[p]}(\mathbb{R}^d \oplus \mathbb{R}^n))$ converging to $\mathbf{x} \in C^{0,\omega,p}(G^{[p]}(\mathbb{R}^d \oplus \mathbb{R}^n))$ in the $d_{\omega,p}$-topology,

$$(2.2) \qquad \sup_{\substack{|y_0| \leq R \\ \|V\|_{\mathrm{Lip}(p+\varepsilon)} \leq R}} d_{\omega,p}(I_{y_0,V}(\mathbf{x}_n), I_{y_0,V}(\mathbf{x})) \underset{n \to \infty}{\to} 0.$$

To see this, one first observes that the continuity of integration along of a one-form is uniform over the set of one-forms with Lipschitz norm bounded by a given $R$. Then, the path $I_{y_0,V}(\mathbf{x})$ over small time is the fixed point of a map (integrating along a one-form) which is a contraction. Reading the estimate in [15], one sees that this map is uniformly, over the set $\{|y_0| \leq R; \|V\|_{\mathrm{Lip}(p+\varepsilon)} \leq R\}$, a contraction with parameter strictly less than 1.

The next theorem was also proved in [15], and deals with the continuity of the flow.

THEOREM 2. *If $(y_0^n)_n$ is a $\mathbb{R}^d$-valued sequence converging to $y_0$, then for all $R > 0$,*

$$\sup_{\substack{\|\mathbf{x}\|_{\omega,p} \leq R \\ \|V\|_{\mathrm{Lip}(p+\varepsilon)} \leq R}} d_{\omega,p}(I_{y_0^n,V}(\mathbf{x}), I_{y_0,V}(\mathbf{x})) \underset{n \to \infty}{\to} 0.$$

The next theorem shows that the Itô map it continuous when one varies the vector fields defining the differential equation. It does not seem to have appeared anywhere, despite that its proof does not involve any new ideas.

THEOREM 3.   *Let $V^n = (V_1^n, \ldots, V_d^n)$ be a sequence of $d$ Lip$(p+\varepsilon)$-vector fields on $\mathbb{R}^n$ such that*

$$\lim_{n \to \infty} \max_i \|V_i^n - V_i\|_{\mathrm{Lip}(p+\varepsilon)} = 0.$$

*Then, for all $R > 0$,*

$$\lim_{n \to \infty} \sup_{\substack{\|\mathbf{x}\|_{\omega,p} \leq R \\ |y_0| \leq R}} d_{\omega,p}(I_{y_0, V^n}(\mathbf{x}), I_{y_0, V}(\mathbf{x})) = 0.$$

PROOF.   We use the notations of [14]. First consider the Lip$(p+\varepsilon-1)$-one-forms $\theta_i : \mathbb{R}^d \to \mathrm{Hom}(\mathbb{R}^d, \mathbb{R}^n)$, $i \in \mathbb{N} \cup \{\infty\}$. We assume that $\theta_n$ converges to $\theta_\infty$ in the $(p+\varepsilon-1)$-Lipschitz topology when $n \to \infty$. For $n \in \mathbb{N} \cup \{\infty\}$, $\int \theta_n(\mathbf{x}) \, d\mathbf{x}$ is the unique rough path associated to the almost multiplicative functional

$$(Z(\mathbf{x})_{s,t}^n)^i = \sum_{l_1, \ldots, l_i = 0}^{[p]-1} d^{l_1}\theta_n(\mathbf{x}_s^1) \otimes \cdots \otimes d^{l_i}\theta_n(\mathbf{x}_s^1) \left( \sum_{\pi \in \Pi_{(l_1, \ldots, l_i)}} \pi(\mathbf{x}_{s,t}^{i+l_1+\cdots+l_i}) \right)$$

(see [14] for the definition of $\Pi$). It is obvious that

$$\sup_{\|\mathbf{x}\|_{\omega,p} \leq R} \max_i \frac{|(Z(\mathbf{x})_{s,t}^n)^i - (Z(\mathbf{x})_{s,t}^\infty)^i|}{\omega(s,t)^{i/p}} \underset{n \to \infty}{\to} 0,$$

which implies by Theorem 3.1.2 in [14] that

$$\sup_{\|\mathbf{x}\|_{\omega,p} \leq R} d_{\omega,p}\left( \int \theta_n(\mathbf{x}) \, d\mathbf{x}, \int \theta_\infty(\mathbf{x}) \, d\mathbf{x} \right) \underset{n \to \infty}{\to} 0.$$

Now consider the Picard iteration sequence $(\mathbf{z}_m^n)_{m \geq 0}$ introduced in [14], formula (4.10), to construct to $I_{y_0, V^n}(\mathbf{x})$:

$$\mathbf{z}_0^n = (0, y_0)$$

$$\mathbf{z}_{m+1}^n = \int h_n(\mathbf{z}_m^n) \, d\mathbf{z}_m^n,$$

where $h_n$ is the one-form defined by the formula $h_n(x,y)(dX, dY) = (dX, V^n(y) \, dX)$ (by $V(y) \, dX$ we mean $\sum V_i(y)^n \, dX^i$). Lyons proved that

$$\lim_{m \to \infty} \sup_{\substack{\|\mathbf{x}\|_{\omega,p} \leq R \\ |y_0| \leq R}} \sup_n d_{\omega,p}(\mathbf{z}_m^n, \mathbf{z}_\infty^n) = 0;$$

we have the supremum over all $n$ here because the $p+\varepsilon$-Lipschitz norm of the $V^n$ are uniformly bounded in $n$. Moreover, we have just seen that for all fixed $m$, $\lim_{n \to \infty} \sup_{|y_0| \leq R} d_{\omega,p}(\mathbf{z}_m^n, \mathbf{z}_m^\infty) = 0$. Therefore, with a $3\varepsilon$-type argument, we obtain our theorem.   $\square$

The three previous theorems actually give that the map

$$(y_0, V, \mathbf{x}) \to I_{y_0, V}(\mathbf{x})$$

is continuous in the product topology $\mathbb{R}^n \times (\text{Lip}(p + \varepsilon)$ on $\mathbb{R}^n)^d \times C^{0,\omega,p}(G^{[p]}(\mathbb{R}^d))$.

In the reminder of this section, $p$ is a real in $[2, 3)$. Now consider a $\text{Lip}(1 + \varepsilon)$-vector field on $\mathbb{R}^n$ denoted $V_0$, and for a path $x$ of bounded variation, we consider $y$ to be the solution of

$$(2.3) \qquad \begin{cases} dy_t = V_0(y_t)\, dt + V(y_t)\, dx_t, \\ y_t = y_0, \end{cases}$$

Lifting $x$ and $(x \oplus y)$ to $S(x)$ and $S(x \oplus y)$ we consider the map which at $S(x)$ associates $S(x, y)$. We denote it $I_{y_0,(V_0,V)}$. The following extension of the universal limit theorem was obtained in [12].

THEOREM 4. *The map $I_{y_0,(V_0,V)}$ is continuous from $(C^{0,\omega,p}(G^{[p]}(\mathbb{R}^d)), d_{\omega,p})$ into $(C^{0,\omega,p}(G^{[p]}(\mathbb{R}^d \oplus \mathbb{R}^n)), d_{\omega,p})$. More precisely, for all $R > 0$,*

$$(2.4) \qquad \sup_{\substack{|y_0| \le R \\ \|V\|_{\text{Lip}(p+\varepsilon)} \le R}} d_{\omega,p}(I_{y_0,(V_0,V)}(\mathbf{x}_n), I_{y_0,(V_0,V)}(\mathbf{x})) \underset{n \to \infty}{\to} 0.$$

Let $x_n$ be a sequence of paths of bounded variation such that $S(x_n)$ converges in the $d_{\omega,p}$-topology to a geometric $p$-rough path $\mathbf{x}$, and define $y_n$ to be the solution of equation (2.3) replacing $x$ by $x_n$. Then, the previous theorem says that $S(x_n \oplus y_n)$ converges in the $d_{\omega,p}$-topology to a geometric $p$-rough path $\mathbf{z}$. We say that $\mathbf{y}$, the projection of $\mathbf{z}$ onto $G^{[p]}(\mathbb{R}^n)$ is the solution of the rough differential equation

$$d\mathbf{y}_t = V_0(\mathbf{y}_t)\, dt + V(\mathbf{y}_t)\, d\mathbf{x}_t$$

with initial condition $y_0$.

We obtain, as before, the following two theorems:

THEOREM 5. *If $(y_0^n)_n$ is a $\mathbb{R}^d$-valued sequence converging to $y_0$, then*

$$\sup_{\substack{\|\mathbf{x}\|_{\omega,p} \le R \\ \|V\|_{\text{Lip}(p+\varepsilon)} \le R \\ \|V_0\|_{\text{Lip}(1+\varepsilon)} \le R}} d_{\omega,p}(I_{y_0^n,(V_0,V)}(\mathbf{x}), I_{y_0,(V_0,V)}(\mathbf{x})) \to_{n \to \infty} 0.$$

THEOREM 6. *Let $(V^n = (V_1^n, \ldots, V_d^n))_{n \ge 0}$ be a sequence of $d$ $\text{Lip}(p+\varepsilon)$-vector fields on $\mathbb{R}^n$ and $(V_0^n)_{n \ge 0}$ a sequence of $\text{Lip}(1 + \varepsilon)$-vector fields on $\mathbb{R}^n$, such that*

$$\lim_{n \to \infty} \max\left\{ \|V_0^n - V_0\|_{\text{Lip}(1+\varepsilon)}, \max_{1 \le i \le d} \|V_i^n - V_i\|_{\text{Lip}(p+\varepsilon)} \right\} = 0.$$

*Then, if* $\mathbf{x} \in C^{0,\omega,p}(G^{[p]}(\mathbb{R}^d))$,

$$\lim_{n \to \infty} \sup_{\substack{\|\mathbf{x}\|_{\omega,p} \le R \\ |y_0| \le R}} d_{\omega,p}(I_{y_0,(V_0^n,V^n)}(\mathbf{x}), I_{y_0,(V_0,V)}(\mathbf{x})) = 0.$$

2.4. *Solving anticipative stochastic differential equations via rough paths.*
We fix a $p \in (2,3)$ and, for simplicity, the control $\omega(s,t) = t - s$ (i.e. we deal with Hölder topologies), although we could have been more general and have considered a wide class of controls as in [6] (i.e., we could have consider modulus type topologies). We define $\mathbf{B}$ the Stratonovich lift to a geometric $p$-rough path of the Brownian motion $B$ with the formula

$$\mathbf{B}_t = \left(1, B_t, \int_0^t B_u \otimes \circ dB_u\right).$$

$\mathbf{B}$ is a $G^2(\mathbb{R}^d)$-valued path, and almost surely, $\|\mathbf{B}\|_{\omega,p} < \infty$.

Consider $V_0$ a random vector field on $\mathbb{R}^n$ almost surely in $\mathrm{Lip}(1 + \varepsilon)$, that is, a measurable map from

$$V_0 : \Omega \times \mathbb{R}^n \to \mathbb{R}^n$$

such that $V_0(\omega, \cdot) \in \mathrm{Lip}(1 + \varepsilon)$ for $\omega$ in a set a full measure, and $V = (V_1, \ldots, V_d)$, where $V_1, \ldots, V_d$ are random vector fields on $\mathbb{R}^n$ almost surely in $\mathrm{Lip}(2 + \varepsilon)$, and a random variable $y_0 \in \mathbb{R}^n$ finite almost surely.

$I_{y_0,(V_0,V)}(\mathbf{B})$ is then almost surely well defined, and its projection onto $G^2(\mathbb{R}^n)$ is the solution, in the rough path sense, of the anticipative stochastic differential equation

$$(2.5) \qquad d\mathbf{y}_t = V_0(\mathbf{y}_t)\,dt + V(\mathbf{y}_t)\,d\mathbf{B}_t$$

with initial condition $y_0$.

The next section introduces the notion of good rough paths sequence, and its properties. Showing that linear approximations of Brownian motion form good rough path sequences (in some sense that will be precise later on) will prove that $\mathbf{y}^1$ is solution of the anticipative Stratonovich stochastic differential equation (1.1). In particular, the solution that we construct coincides with the one constructed in [20].

## 3. Good rough path sequence.

3.1. *Definitions.* We fix a parameter $p > 2$, and a control $\omega$. $\mathbb{R}^d$ and $\widetilde{\mathbb{R}^d}$ will denote two identical copies of $\mathbb{R}^d$.

Let $p > 2$, and $q$ such that $1/p + 1/q > 1$. We consider $x$ and $y$ two $\mathbb{R}^d$-valued paths of bounded variation. We let $\mathbf{y} = S(y)$ to be the canonical lift of $y$ to a $G^{[p]}(\mathbb{R}^d)$-valued path. We let

$$(3.1) \qquad S'(x, S(y)) := S(x \oplus y)$$

be the canonical lift of $x \oplus y$ to a $G^{[p]}(\mathbb{R}^d \oplus \mathbb{R}^d)$-valued path and

(3.2) $$S''(S(x)) := S'(x, S(x)) = S(x \oplus x)$$

be the canonical lift of $x \oplus x$ to a $G^{[p]}(\mathbb{R}^d \oplus \mathbb{R}^d)$-valued path.

PROPOSITION 1. *Let $x$ be a $\mathbb{R}^d$-valued path of finite $q$-variation, and $\mathbf{y}$ a $G^{[p]}(\mathbb{R}^d)$-valued path of finite $p$-variation. Let $(x_n, y_n)$ be a sequence of $\mathbb{R}^d \oplus \mathbb{R}^d$-valued path such that $d_{\omega,p}(x_n, x) \to_{n \to \infty} 0$ and $d_{\omega,p}(S(y_n), \mathbf{y}) \to_{n \to \infty} 0$. Then:*

*(i) $S'(x_n, S(y_n))$ converges in the $d_{\omega,p}$-topology, and the limit is independent of the choice of the sequence $(x_n, y_n)$. We denote this limit element $S'(x, \mathbf{y})$.*

*(ii) $S''(S(y_n))$ converges in $d_{\omega,p}$-topology, and the limit is independent of the choice of the sequence $(y_n)$. We denote this element $S''(\mathbf{y}, \mathbf{y})$.*

PROOF. This is simply obtained using Theorem 3.1.2 in [14], which says that the procedure which at an almost multiplicative functional associates a rough path is continuous, and we leave the details to the reader. □

EXAMPLE 2. If $2 \le p < 3$,

$$S'(x, \mathbf{y})_t = \left(1, x_t \oplus \mathbf{y}_t^1, \int_0^t x_u \otimes dx_u \oplus \int_0^t x_u \otimes d\mathbf{y}_u^1 \oplus \int_0^t \mathbf{y}_u^1 \otimes dx_u \oplus \mathbf{y}_t^2\right).$$

The three integrals are well defined Young's integrals [26].

We introduce the notion of a good $p$-rough path sequence.

DEFINITION 2. Let $(x_n)_n$ be a sequence of $\mathbb{R}^d$-valued paths of bounded variation, and $\mathbf{x}$ a geometric $p$-rough path. We say that $(x_n)_{n \in \mathbb{N}}$ is a good $p$-rough path sequence (associated to $\mathbf{x}$) (for the control $\omega$) if and only if

(3.3) $$\lim_{n \to \infty} d_{\omega,p}(S'(x_n, \mathbf{x}), S''(\mathbf{x})) = 0.$$

In particular, if $(x_n)_n$ is a good $p$-rough path sequence associated to $\mathbf{x}$, for the control $\omega$, $x_n$ converges to $\mathbf{x}$ in the topology induced by $d_{\omega,p}$.

PROPOSITION 2. *Assume $2 \le p < 3$. The sequence $(x(n))_n$ of paths of bounded variation is good rough path sequence associated to $\mathbf{x}$, for the control $\omega$, if and only if*

$$\lim_{n \to \infty} \sup_{0 \le s < t \le 1} \frac{|x(n)_{s,t} - \mathbf{x}_{s,t}^1|}{\omega(s,t)^{1/p}} = 0,$$

$$\lim_{n \to \infty} \sup_{0 \le s < t \le 1} \frac{|\int_s^t x(n)_{s,u} \otimes dx(n)_u - \mathbf{x}_{s,t}^2|}{\omega(s,t)^{2/p}} = 0,$$

$$\lim_{n \to \infty} \sup_{0 \le s < t \le 1} \frac{|\int_s^t \mathbf{x}_{s,u}^1 \otimes dx(n)_u - \mathbf{x}_{s,t}^2|}{\omega(s,t)^{2/p}} = 0.$$

PROOF.   $d_{\omega,p}(S'(x(n), \mathbf{x}), S''(\mathbf{x}, \mathbf{x}))$ goes to zero if and only if

$$\max\{A_1^n, \sqrt{A_2^n}, \sqrt{A_3^n}, \sqrt{A_4^n}\}$$

goes to zero, where

$$A_1^n = \sup_{0 \le s < t \le 1} \frac{|\mathbf{x}_{s,t}^1 - x(n)_{s,t}|}{\omega(s,t)^{1/p}},$$

$$A_2^n = \sup_{0 \le s < t \le 1} \frac{|\int_s^t x(n)_{s,u} \otimes dx(n)_u - \mathbf{x}_{s,t}^2|}{\omega(s,t)^{2/p}},$$

$$A_3^n = \sup_{0 \le s < t \le 1} \frac{|\int_s^t x(n)_{s,u} \otimes d\mathbf{x}_u^1 - \mathbf{x}_{s,t}^2|}{\omega(s,t)^{2/p}},$$

$$A_4^n = \sup_{0 \le s < t \le 1} \frac{|\int_s^t \mathbf{x}_{s,u}^1 \otimes dx(n)_u - \mathbf{x}_{s,t}^2|}{\omega(s,t)^{2/p}}.$$

By integration by part, we see that

$$\int_s^t \mathbf{x}_{s,u}^1 \otimes dx(n)_u - \mathbf{x}_{s,t}^2 = \mathbf{x}_{s,t}^1 \otimes x(n)_{s,t} - \int_s^t d\mathbf{x}_{s,u}^1 \otimes x(n)_u - \mathbf{x}_{s,t}^2$$

$$= \mathbf{x}_{s,t}^1 \otimes (x(n)_{s,t} - \mathbf{x}_{s,t}^1) - \int_s^t d\mathbf{x}_{s,u}^1 \otimes x(n)_u$$

$$+ ((\mathbf{x}_{s,t}^1)^{\otimes 2} - \mathbf{x}_{s,t}^2)$$

Let $\pi$ be the linear operator from $\mathbb{R}^d \otimes \mathbb{R}^d$ onto itself defined by $\pi(x \otimes y) = \pi(y \otimes x)$. Observe that $\pi(\mathbf{x}_{s,t}^2) = (\mathbf{x}_{s,t}^1)^{\otimes 2} - \mathbf{x}_{s,t}^2$. Also, we have for all $z \in \mathbb{R}^d \otimes \mathbb{R}^d$, $|z| = |\pi(z)|$. Hence, we obtain from this remark and the above equation that

$$\int_s^t x(n)_{s,u} \otimes d\mathbf{x}_u^1 - \mathbf{x}_{s,t}^2 = \pi\left(\int_s^t \mathbf{x}_{s,u}^1 \otimes dx(n)_u - \mathbf{x}_{s,t}^2\right)$$

$$+ (x(n)_{s,t} - \mathbf{x}_{s,t}^1) \otimes \mathbf{x}_{s,t}^1$$

This implies that $A_3^n \le \|\mathbf{x}\|_{\omega,p} \cdot (A_1^n + A_4^n)$, which proves the proposition.   $\square$

By definition of a geometric $p$-rough path, there always exists a sequence of smooth $x_n$ such that $S(x_n)$ converges to $\mathbf{x}$ in the topology induced by $d_{p\text{-var}}$. However, this does not imply that equation (3.3) holds (or equivalently for $[p] = 2$: that the conditions in Proposition 2 hold). We now give an example of a geometric rough path $\mathbf{x}$ for which there exists no good sequence associated to it.

Example 3.   Consider

$$\mathbf{x}_t = \exp(t[e_1, e_2])$$
$$= \exp(t(e_1 \otimes e_2 - e_2 \otimes e_1))$$

where $e_1, e_2$ is a basis of $\mathbb{R}^2$. Note that $\mathbf{x}$ is a geometric $p$-rough path and

$$\mathbf{x}_t^1 \equiv 0 \in \mathbb{R}^2, \qquad x_t^2 \equiv t[e_1, e_2] \in \mathbb{R}^2 \otimes \mathbb{R}^2.$$

If $(x_n)$ is a sequence of smooth paths in $\mathbb{R}^2$, then $\int_0^t \mathbf{x}^1 \otimes dx_n \equiv 0$ trivially converges to $0 \in \mathbb{R}^2 \otimes \mathbb{R}^2$. On the other hand, if $x_n$ was a good $p$-rough path sequence associated to $\mathbf{x}$ it should converge to $t[e_1, e_2]$

3.2. *Stability of good rough path sequences.*   Nonetheless, there are some good news. Good rough path sequences are stable under integration, as shows the following theorem. We fix a control $\omega$.

Theorem 7.   *Let $(x_n)_{n \in \mathbb{N}}$ denote a good $p$-rough path sequence associated to $\mathbf{x}$ for the control $\omega$, and $\theta : \mathbb{R}^d \to \mathrm{Hom}(\mathbb{R}^d, \mathbb{R}^n)$ be a $\mathrm{Lip}(p - 1 + \varepsilon)$ one-form, with $\varepsilon > 0$. Then $(\int \theta(x_n) \, dx_n)_n$ is a good rough path sequence associated to $\int \theta(\mathbf{x}) \, d\mathbf{x}$ for the control $\omega$. Moreover,*

$$(3.4) \qquad \lim_{n \to \infty} S\left( \int \theta(\mathbf{x}^1) \, dx_n \right) = \int \theta(\mathbf{x}) \, d\mathbf{x}.$$

Proof.   Consider the $\mathrm{Lip}(p - 1 + \varepsilon)$ one form

$$\widetilde{\theta} : \mathbb{R}^d \oplus \mathbb{R}^d \to \mathrm{Hom}(\mathbb{R}^d \oplus \mathbb{R}^d, \mathbb{R}^n \oplus \mathbb{R}^n),$$
$$((x, y), (dX, dY)) \to \theta(x) \, dX \oplus \theta(y) \, dY.$$

All the limits considered in this proof are to be understood to be in the topology induced by $d_{\omega, p}$. By the continuity of the integral [14, 15] and of the operator $S'$ (by Proposition 1),

$$S'\left( \int \theta(x_n) \, dx_n, \int \theta(\mathbf{x}) \, d\mathbf{x} \right) = \lim_{m \to \infty} S'\left( \int \theta(x_n) \, dx_n, S\left( \int \theta(x_m) \, dx_m \right) \right)$$
$$= \lim_{m \to \infty} S\left( \int \theta(x_n) \, dx_n \oplus \int \theta(x_m) \, dx_m \right)$$
$$= \lim_{m \to \infty} S\left( \int \widetilde{\theta}(x_n, x_m) \, d(x_n, x_m) \right)$$

$$= \lim_{m \to \infty} \int \tilde{\theta}(S(x_n \oplus x_m)) \, dS(x_n \oplus x_m)$$

$$= \int \widetilde{\theta}(S'(x_n, \mathbf{x})) \, dS'(x_n, \mathbf{x}).$$

In in the last line we have used once again the continuity of the integral and the assumption that $(x_n)$ is a good $p$-rough path sequence associated to $\mathbf{x}$.

Therefore, $S'(\int \theta(x_n) \, dx_n, \int \theta(\mathbf{x}) \, d\mathbf{x})$ converges when $n \to \infty$ to

$$\int \widetilde{\theta}(S''(\mathbf{x})) \, dS''(\mathbf{x}) = \lim_n \int \widetilde{\theta}(S''(S(x_n))) \, dS''(S(x_n))$$

$$= \lim_n \int \widetilde{\theta}(S(x_n \oplus x_n)) \, dS(x_n \oplus x_n)$$

$$= \lim_n S\left( \int \theta(x_n) \, dx_n \oplus \int \theta(x_n) \, dx_n \right)$$

$$= S''\left( \int \theta(\mathbf{x}) \, d\mathbf{x} \right),$$

which proves the first assertion.

For the second one, consider the map

$$\widehat{\theta} \colon \mathbb{R}^d \oplus \mathbb{R}^d \to \mathrm{Hom}(\mathbb{R}^d \oplus \mathbb{R}^d, \mathbb{R}^n),$$

$$((x, y), (dX, dY)) \to \theta(y) \, dX.$$

By the continuity of the integral, we obtain that $\int \widehat{\theta}(S'(x_n, \mathbf{x})) \, dS'(x_n, \mathbf{x})$ converges in the topology induced by $d_{\omega, p}$ to $\int \widehat{\theta}(S''(\mathbf{x})) \, dS''(\mathbf{x})$. This is our result as $\int \widehat{\theta}(S'(x_n, \mathbf{x})) \, dS'(x_n, \mathbf{x}) = S(\int \theta(\mathbf{x}) \, dx_n)$. $\quad\square$

We show now that good rough path sequences are stable under the Itô map.

THEOREM 8. *Let $V_1, \ldots, V_d$ be $d$ elements of* $\mathrm{Lip}(p + \varepsilon)$-*vector fields on* $\mathbb{R}^n$, *and $V = (V_1, \ldots, V_d)$ thought as a linear map from $\mathbb{R}^d$ into* $\mathrm{Lip}(p + \varepsilon)$ *vector fields on $\mathbb{R}^n$. Assume that $(x_n)_n$ is a good $p$-rough path sequence associated to $\mathbf{x}$ for the control $\omega$. Denote by $y_n$ the solution of the ordinary differential equation*

$$\begin{cases} dy_n(t) = V(y_n(t)) \, dx_n(t), \\ y_n(0) = y_0. \end{cases}$$

*Then $(x_n \oplus y_n)_n$ is a good $p$-rough path sequence associated to $\mathbf{z} = I_{y_0, V}(\mathbf{x})$ for the control $\omega$.*

Proof. The proof is similar to the proof of Theorem 7. Denote by $\widetilde{V}$ the linear map from $\mathbb{R}^d \oplus \mathbb{R}^d$ into the $\mathrm{Lip}(p+\varepsilon)$ vector fields on $\mathbb{R}^n$ by the formula $\widetilde{V}(y_1, y_2)(dx_1, dx_2) = (V(y_1)\, dx_1, V(y_2)\, dx_2)$. First notice that

$$I_{(y_0,y_0),\widetilde{V}}(S''(\mathbf{x})) = S''(I_{y_0,V}(\mathbf{x}))$$

and that

$$S'(z_n, \mathbf{z}) = I_{(y_0,y_0),\widetilde{V}}(S'(x_n, \mathbf{x})).$$

Hence, since $(x_n)$ is a good sequence and by continuity of the Itô map, $S'(z_n, \mathbf{z})$ converges as $n \to \infty$ to

$$I_{(y_0,y_0),\widetilde{V}}(S''(\mathbf{x})) = S''(I_{y_0,V}(\mathbf{x})) = S''(\mathbf{z}). \qquad \square$$

From our two previous theorems, we immediately obtain the following corollary:

Corollary 1. *We keep the notation of Theorem 8. Then, in the topology induced by* $d_{\omega,p}$,

$$\mathbf{y} = \lim_{n\to\infty} S\left(y_0 + \int V(\mathbf{y}_u^1)\, dx_n(u)\right).$$

*In particular, looking at the first level of this equation, we obtain that*

$$y_0 + \int_0^{\cdot} V(y_u)\, dx_n(u) \to y.$$

Remark 1. In the previous theorem and its corollary, with no modification in the proof, for $2 \le p < 3$, one can obtain the same results replacing the map $I_{y_0,V}$ by the map $I_{y_0,(V_0,V)}$, where $V_0$ is a $\mathrm{Lip}(1+\varepsilon)$-vector fields on $\mathbb{R}^n$. In other words, one can consider differential equations with a time drift and (almost) minimal smoothness condition on $V_0$.

3.3. *Piecewise-linear approximation of Brownian motion as a good rough path sequence.* We fix a $p \in (2,3)$ and $\omega(s,t) = t - s$. We recall from the introduction that $B$ is a $d$-dimensional Brownian motion, and that $\mathbf{B}$ is the Stratonovich lift of $B$ to a geometric $p$-rough path. Let $B^D$ be the $D$-linear approximation of $B$ [equation (1.2)]. Let $(D^n)$ be a sequence of subdivisions of $[0,1]$ which steps tends to 0 when $n$ tends to infinity. If $D^n = (\frac{k}{2^n}, 0 \le k \le 2^n)$, we know from [15] that, almost surely, $S(B^{D^n})$ converges in $p$-variation to $\mathbf{B}$. If $D^n$ is an increasing sequence of subdivision, that is, if $D^{n+1} \subset D^n$ for all $n$, a martingale argument proved in [6] that almost surely and in $L^q$, $q \ge 1$, $S(B^{D^n})$ converges in $1/p$-Hölder distance (and even some modulus distances) to $\mathbf{B}$. The following theorem goes a bit deeper in the convergence of piecewise linear approximations of the Brownian motion.

THEOREM 9.  *Let $(D^n)$ be a sequence of subdivision which steps size tends to 0. Then $d_{\omega,p}(S'(B^{D^n}, \mathbf{B}), S''(\mathbf{B}, \mathbf{B}))$ converges when $n$ tends to infinity to 0 in $L^q$, $q \geq 1$ and in probability.*

*If $D^n = (\frac{k}{2^n}, 0 \leq k \leq 2^n)$, the convergence also holds almost surely, that is, $B^n = B^{D^n}$ is almost surely a good p-rough path sequence associated to $\mathbf{B}$.*

We decompose the proof in four lemmas.

LEMMA 1.  *For all $q > p > 1$, the $L^q$ and $L^p$-norms on the $k^{th}$ Wiener chaos are equivalent.*

PROOF.  This is a simple consequence of the hypercontractivity of the Ornstein–Uhlenbeck semigroup; see [19], page 57, for example.  □

LEMMA 2.  *Let $D$ be a subdivision of $[0,1]$. Then, for all $q \geq 1$ and $p' > 2$, there exists $\mu > 0$ and $C_{p',q,\mu} < \infty$ such that for all $s < t \in D$*

$$\left\| \int_s^t B_{s,u} \otimes dB_u^D - \mathbf{B}_{s,t}^2 \right\|_{L^q} \leq C_{p',q,\mu} |D|^\mu |t-s|^{2/p'},$$

$$\left\| \int_s^t B_{s,u}^D \otimes dB_{s,u}^D - \mathbf{B}_{s,t}^2 \right\|_{L^q} \leq C_{p',q,\mu} |D|^\mu |t-s|^{2/p'}$$

PROOF.  From the previous lemma, we can take $q = 2$. We write $(s,t) = (t_m, t_n)$, with $0 \leq m < n \leq |D|$, where $D = (t_i)_{0 \leq i \leq D}$. It is easy to see that

$$\int_s^t B_{s,u} \otimes dB_u^D - \mathbf{B}_{s,t}^2 = \sum_{k=m}^{n-1} \left( \int_{t_k}^{t_{k+1}} B_{t_k,u} \otimes dB_u^D - \mathbf{B}_{t_k,t_{k+1}}^2 \right),$$

$$\int_s^t B_{s,u}^D \otimes dB_{s,u}^D - \mathbf{B}_{s,t}^2 = \sum_{k=m}^{n-1} ((B_{t_k,t_{k+1}})^{\otimes 2} - \mathbf{B}_{t_k,t_{k+1}}^2).$$

Therefore, by independence of increment,

$$\left\| \int_s^t B_{s,u} \otimes dB_u^D - \mathbf{B}_{s,t}^2 \right\|_{L^2}^2 = \sum_{k=m}^{n-1} E\left( \left( \int_{t_k}^{t_{k+1}} B_{t_k,u} \otimes dB_u^D - \mathbf{B}_{t_k,t_{k+1}}^2 \right)^2 \right)$$

$$= C \sum_{k=m}^{n-1} (t_{k+1} - t_k)^2$$

$$\leq C|D|^{2-4/p'} \sum_{k=m}^{n-1} (t_{k+1} - t_k)^{4/p'}$$

$$\leq C|D|^{2-4/p'} (t-s)^{4/p'}.$$

We also obtain the same estimate for $\| \int_s^t B_{s,u}^D \otimes dB_{s,u}^D - \mathbf{B}_{s,t}^2 \|_{L^q}$, which concludes the proof. $\square$

LEMMA 3. *Let $D$ be a subdivision of $[0,1]$. Then, for all $q \geq 1$, there exists $\mu > 0$ and $C_{q,\mu} < \infty$ such that*

$$(3.5) \qquad \left\| \sup_{0 \leq s < t \leq 1} \frac{|\int_s^t B_{s,u} \otimes dB_u^D - \mathbf{B}_{s,t}^2|}{|t-s|^{2/p}} \right\|_{L^q} \leq C_{q,\mu} |D|^\mu,$$

$$(3.6) \qquad \left\| \sup_{0 \leq s < t \leq 1} \frac{|\int_s^t B_{s,u}^D \otimes dB_{s,u}^D - \mathbf{B}_{s,t}^2|}{|t-s|^{2/p}} \right\|_{L^q} \leq C_{q,\mu} |D|^\mu.$$

PROOF. We only prove equation (3.5) as the proof for the other estimates is similar. We define $X_{s,t}^D = \int_s^t B_{s,u} \otimes dB_u^D - \mathbf{B}_{s,t}^2$, where $D = (0 = t_0 \leq t_1 < \cdots < t_{|D|^*} = 1)$ is fixed subdivision of $[0,1]$. First assume that $t_i \leq s < t \leq t_{i+1}$. We let $p' \in (2,p)$; it is easy to check that there exists $C < \infty$ independent of $D$ such that $\|B^D\|_{\omega,p'} \leq C\|B\|_{\omega,p'}$. Therefore,

$$\frac{|X_{s,t}^D|}{|t-s|^{2/p}} = \frac{|\mathbf{B}_{s,t}^2| + |(\int_s^t B_{s,u}\,du) \otimes B_{t_i,t_{i+1}}/(t_{i+1}-t_i)|}{|t-s|^{2/p}}$$

$$\leq \frac{C\|B\|_{\omega,p'}^2(|t-s|^{2/p'} + (\int_s^t (u-s)^{1/p'}\,du)(t_{i+1}-t_i)^{1/p'-1})}{|t-s|^{2/p}}$$

$$\leq C\|B\|_{\omega,p'}^2(|t-s|^{2/p'-2/p} + (t-s)^{1+1/p'-2/p}(t_{i+1}-t_i)^{1/p'-1}).$$

Bounding $t-s$ and $t_{i+1} - t_i$ by $|D|$, the mesh size of $D$, we obtain that

$$(3.7) \qquad \max_{k \in \{0,\ldots,|D|^*\}} \sup_{t_k \leq s < t \leq t_{k+1}} \frac{|X_{s,t}^D|}{|t-s|^{2/p}} \leq C\|B\|_{\omega,p'}^2 |D|^{2(1/p'-1/p)}.$$

Defining $t_D$ to be the biggest real in $D$ less than or equal to $t$, and $s^D$ the smallest real in $D$ greater than to $s$, the above estimate rewrites

$$(3.8) \qquad \sup_{\substack{0 \leq s < t \leq 1 \\ t_D < s^D}} \frac{|X_{s,t}^D|}{|t-s|^{2/p}} \leq C\|B\|_{\omega,p'}^2 |D|^{2(1/p'-1/p)}$$

and an $L^q$-estimate is immediate. Hence we are left to prove that

$$\left\| \sup_{\substack{0 \leq s < t \leq 1 \\ s^D \leq t_D}} \frac{|X_{s,t}^D|}{|t-s|^{2/p}} \right\|_{L^q} \leq C_\mu |D|^\mu.$$

Now observe that for all $s < t < u$,

$$(3.9) \qquad X_{s,u}^D = X_{s,t}^D + X_{t,u}^D + B_{s,t} \otimes (B_{t,u}^D - B_{t,u}).$$

Hence, for all $s < t$ such that $s^D \leq t_D$, so that

$$X_{s,t}^D = X_{s,s^D}^D + X_{s^D,t_D}^D + X_{t_D,t}^D + B_{s,t_D} \otimes (B_{t_D,t}^D - B_{t_D,t}).$$

From (3.8)

$$\left\| \sup_{\substack{0 \leq s < t \leq 1 \\ s^D \leq t_D}} \frac{|X_{s,s^D}^D|}{|s^D - s|^{2/p}} + \frac{|X_{t_D,t}^D|}{|t - t_D|^{2/p}} \right\|_{L^q} \leq C_\mu |D|^\mu.$$

Compatibility of tensor norms shows—similar as above but easier—

$$\frac{|B_{s,t_D} \otimes (B_{t_D,t}^D - B_{t_D,t})|}{|t - s|^{2/p}} \leq C \|B\|_{\omega,p'}^2 |D|^{1/p' - 1/p};$$

therefore, we just need to check that for some $\mu > 0$,

(3.10) $\qquad \left\| \sup_{0 \leq s < t \leq 1} \frac{|X_{s^D,t_D}^D|}{|t - s|^{2/p}} \right\|_{L^q} \overset{\text{trivial}}{\leq} \left\| \max_{s < t \in D} \frac{|X_{s,t}^D|}{|t - s|^{2/p}} \right\|_{L^q} \leq C_\mu |D|^\mu.$

Now consider $\widetilde{X}^D : [0,1] \to \mathbb{R}^d \otimes \mathbb{R}^d$ the linear path on the intervals $[t_i, t_{i+1}]$, such that for all $i$, $\widetilde{X}_{t_{i+1}}^D - \widetilde{X}_{t_i}^D = X_{t_i,t_{i+1}}^D$, that is, for all $t$,

$$\widetilde{X}_t^D = \widetilde{X}_{t_D}^D + \frac{t - t_D}{t^D - t_D} X_{t_D,t^D}^D.$$

The previous lemma showed that for all $s, t \in D$, $\|\widetilde{X}_t^D - \widetilde{X}_s^D\|_{L^q} \leq C_{p',q,\mu} |D|^\mu \times |t - s|^{2/p'}$. From this, it is easy to check that (changing the constant), this equality remains true for all $s, t \in [0,1]$. Define $B_\alpha = \int_0^1 \int_0^1 |\frac{\widetilde{X}_t^D - \widetilde{X}_s^D}{(t-s)^\alpha}|^{2m} \, ds \, dt$. If $\alpha < \frac{2}{p'} + \frac{1}{2m}$, then $\mathbb{E}(B_\alpha) \leq C_m |D|^{\mu m}$. By Garsia, Rodemich and Rumsey's theorem [24] (used here classical on the normed Euclidean space $\mathbb{R}^d \otimes \mathbb{R}^d$), we obtain that

$$\sup_{0 \leq s < t \leq 1} \frac{|\widetilde{X}_t^D - \widetilde{X}_s^D|}{(t - s)^{\alpha - 1/m}} \leq C B_\alpha^{1/2m}.$$

Therefore, there exists two constants $\vartheta > 0$ and $C_\vartheta < \infty$ (independent of $D$) such that

$$\left\| \max_{s < t \in D} \frac{|X_{s,t}^D|}{|t - s|^{2/p}} \right\|_{L^q} \leq C_\vartheta |D|^\vartheta. \qquad \qquad \square$$

LEMMA 4. *Let $D$ be a subdivision of $[0,1]$. Then, for all $q \geq 1$, there exists $\nu > 0$ and $C_{q,\nu} < \infty$ such that*

$$\|d_{\omega,p}(S'(B^D, \mathbf{B}), S''(\mathbf{B}, \mathbf{B}))\|_{L^q} \leq C_{q,\nu} |D|^\nu.$$

Proof. It is easy to check that $\sup_{0 \le s < t \le 1} \frac{|\mathbf{B}^1_{s,t} - B^D_{s,t}|}{|t-s|^{1/p}} \le C \|\mathbf{B}\|_{\omega,q} |D|^{1/q-1/p}$, where $2 < q < p$. The previous lemmae together with Proposition 2 give the result. $\square$

We can now turn to the proof of Theorem 9.

Proof of Theorem 9. The first part is obvious from the previous lemma. From the second part,

$$\mathbb{P}\left(d_{\omega,p}(S'(B^n, \mathbf{B}), S''(\mathbf{B}, \mathbf{B})) \ge \frac{1}{n}\right) \le n^2 \|d_{\omega,p}(S'(B^n, \mathbf{B}), S''(\mathbf{B}, \mathbf{B}))\|^2_{L^2}$$
$$\le C_{q,\nu} n^2 2^{-n\nu}.$$

Hence, by Borel–Cantelli's lemma, we obtain that, almost surely, $S'(B^n, \mathbf{B})$ converges to $S''(\mathbf{B}, \mathbf{B})$ in the topology induced by $d_{\omega,p}$. $\square$

## 4. Anticipative stochastic analysis.
We present a few applications, such as Wong–Zakai results, support theorems and large deviations in the context of anticipating stochastic calculus.

4.1. *Rough paths solution equals Stratonovich solution.* Fix $p \in (2,3)$ and $\omega(s,t) = t - s$. Consider random vector-fields $V_0$ and $V = (V_1, \ldots, V_d)$, where $V_0$ (resp. $V_1, \ldots, V_d$) is almost surely a Lip$(1 + \varepsilon)$ [resp. Lip$(2 + \varepsilon)$] vector field on $\mathbb{R}^n$, and an a.s. finite random variable $y_0 \in \mathbb{R}^n$. As earlier, $\mathbf{B}$ denotes the canonical lift of Brownian motion to a geometric $p$-rough path. Then, there exists a unique rough path solution of the (anticipative) stochastic differential equation

$$d\mathbf{y}_t = V_0(\mathbf{y}_t) \, dt + V(\mathbf{y}_t) \, d\mathbf{B}_t$$

with (random) initial condition $y_0$.

Theorem 10. *Let $\mathbf{y}^1$ denote the projection to path-level of the rough-path $\mathbf{y}$. Then $\mathbf{y}^1$ solves the anticipative Stratonovich stochastic differential equation,*

$$(4.1) \qquad \mathbf{y}^1_t = y_0 + \int_0^t V_0(\mathbf{y}^1_u) \, du + \int_0^t V(\mathbf{y}^1_u) \circ dB_u.$$

*Moreover, almost surely, $\mathbf{y}^1$ is the limit in $1/p$-Hölder topology of*

$$(4.2) \qquad t \to y_0 + \int_0^t V_0(\mathbf{y}^1_u) \, du + \int_0^t V(\mathbf{y}^1_u) \, dB^n_u,$$

*where $B^n$ is the dyadic linear approximation of $B$ of level $n$.*

PROOF. For any sequence of subdivisions $(D^n)_n$ which mesh size tends to 0 when $n \to \infty$, $d_{\omega,p}(S'(B^{D^n}, \mathbf{B}), S''(\mathbf{B}, \mathbf{B}))$ converges in probability to 0 (Theorem 9) hence by Corollary 1,

$$d_{\omega,p}\left(S\left(y_0 + \int_0^t V_0(\mathbf{y}_u^1) \, du + \int_0^t V(\mathbf{y}_u^1) \, dB_u^{D^n}\right), \mathbf{y}\right)$$

converges in probability to 0. The first level of this equation says precisely that $\mathbf{y}^1$ is solution of the anticipative Stratonovich stochastic differential equation (4.2).

The same argument and the fact that $d_{\omega,p}(S'(B^n, \mathbf{B}), S''(\mathbf{B}, \mathbf{B}))$ converges almost surely to 0 gives the second part of the theorem. □

4.2. *A Wong–Zakai theorem.* The universal limit theorem gives us for free a Wong–Zakai theorem for our solution of the Stratonovich differential equation.

THEOREM 11. *Under the same assumptions and notation than the previous theorem, let $y^n$ be the solution of the differential equation*

$$\begin{cases} dy_t^n = V_0(y_t^n) \, dt + V(y_t^n) \, dB^n(t), \\ y_0^n = y_0. \end{cases}$$

*Then almost surely, $y^n$ converges to $y$ in the topology induced by the $1/p$-Hölder distance.*

PROOF. Let $\mathbf{z}^n = I_{V,y_0}(S(B^n))$. Then by the continuity of the Itô map, $d_{\omega,p}(\mathbf{z}^n, \mathbf{z}) \to 0$. As $\mathbf{z}^n$ projects down onto $y^n$, this is a stronger result that the stated theorem. □

4.3. *A large deviation principle.* Due to the universal limit theorem, Freidlin–Wentzell's type theorems have easy proofs via rough paths [7, 16]. We give here an extension of the Freidlin–Wentzell's theorem, and of the main theorem in [18].

THEOREM 12. *Let $(y_0^\alpha)_{\alpha \geq 0}$ be a family of random elements of $\mathbb{R}^d$, $(V_0^\alpha)_{\alpha \geq 0}$ be a family of random $\mathrm{Lip}(1 + \varepsilon)$-vector fields, $(V_1^\alpha, \ldots, V_d^\alpha)_{\alpha \geq 0}$, $i = 1, \ldots, d$, $d$ families of random $\mathrm{Lip}(2 + \varepsilon)$-vector fields, such that for all $\beta > 0$*

$$\lim_{\alpha \to 0} \alpha \log \mathbb{P}(\|V_0^\alpha - V_0^0\|_{\mathrm{Lip}(1+\varepsilon)} > \beta) = -\infty,$$

$$\lim_{\alpha \to 0} \alpha \log \mathbb{P}\left(\max_{1 \leq i \leq d} \|V_i^\alpha - V_i^0\|_{\mathrm{Lip}(2+\varepsilon)} > \beta\right) = -\infty,$$

$$\lim_{\alpha \to 0} \alpha \log \mathbb{P}(|y_0^\alpha - y_0^0| > \beta) = -\infty.$$

*Define $\mathbf{y}^\alpha$ to be the rough path solution of the differential equation*

$$\begin{cases} d\mathbf{y}_t^\alpha = V_0^\alpha(\mathbf{y}_t^\alpha)\,dt + \sqrt{\alpha}V^\alpha(\mathbf{y}_t^\alpha)\,d\mathbf{B}_t, \\ \mathbf{y}_0^\alpha = y_0^\alpha. \end{cases}$$

*Then $(\mathbf{y}^\alpha)_{\alpha>0}$ satisfies a large deviation principle in the topology induced by $d_{\omega,p}$ with good rate function*

$$J(\mathbf{x}) = \inf_{I_{y_0^0,V^0}(\mathbf{y}) = \mathbf{x}} I(\mathbf{y}),$$

*where*

$$I(\mathbf{x}) = \begin{cases} \frac{1}{2}\int_0^1 |x_u'|^2\,du, & \text{if } S(x) = \mathbf{x} \text{ for some } x \in W^{1,2} \text{ ,} \\ +\infty, & \text{otherwise.} \end{cases}$$

PROOF. In [6], it was proved that $(\delta_{\sqrt{\alpha}}\mathbf{B})_{\alpha>0}$ satisfies a large deviation principle in the topology induced by $d_{\omega,p}$ with good rate function $I$. The assumptions on the vector fields and the initial conditions give that $((V_i^\alpha)_{0\leq i\leq d}, y_0^\alpha, \delta_{\sqrt{\alpha}}\mathbf{B})$ satisfies a large deviation (in the topology product $\mathrm{Lip}(1+\varepsilon), \mathrm{Lip}(2+\varepsilon)^d, |\cdot|$, and $d_{\omega,p}$) with good rate function

$$\begin{cases} \frac{1}{2}\int_0^1 |x_u'|^2\,du, & \text{if } S(x) = \mathbf{x} \text{ for some } x \in W^{1,2},\ V_i = V_i^0, \\ & 1 \leq i \leq d \text{ and } y_0 = y_0^0, \\ +\infty, & \text{otherwise.} \end{cases}$$

By the continuity of the Itô map (Theorems 4, 5 and 6), we obtain our large deviation principle. □

4.4. *Support theorem.* We recall the support theorem for the enhanced Brownian motion (see [4, 5, 6, 16]).

THEOREM 13. *Let and $p > 2$. The support of the law of $\mathbf{B}$ in the $d_{\omega,p}$-topology is the set of paths starting at 0 in $C^{0,\omega,p}(G^{[p]}(\mathbb{R}^d))$ where $\omega(s,t) = t - s$.*

Denote by $I$ the map which maps $(y_0, (V_0, V), \mathbf{x})$ to $I_{y_0,(V_0,V)}(\mathbf{x})$. The following proposition is an obvious corollary of the continuity of the map $I$ (Theorems 4, 5, 6) and of Theorem 13.

PROPOSITION 3. *Let $\mathbf{y}$ be the solution of the rough differential equation*

$$\begin{cases} d\mathbf{y}_t = V_0(\mathbf{y}_t)\,dt + V(\mathbf{y}_t)\,d\mathbf{B}(t), \\ \mathbf{y}_0^1 = y_0, \end{cases}$$

*where $V_0$ is almost surely a* Lip$(1 + \varepsilon)$ *vector field and $V_i$, $i \in \{1, \ldots, d\}$ are almost surely* Lip$(2 + \varepsilon)$ *vector fields, $y_0 \in \mathbb{R}^d$ is almost surely finite. The support of the law of $\mathbf{y}$ in the $d_{\omega,p}$-topology is equal to the image by the map $I$ of the support of the law of $(y_0, V, \mathbf{B})$, in the product of the Euclidean, Lipschitz and $d_{\omega,p}$ topology.*

In particular, if $y_0$ and the vector fields $V_i$ are deterministic, in the $d_{\omega,p}$-topology, the support of the law of $\mathbf{y}$ is equal to the set $I_{y_0,(V_0,V)}(C^{0,\omega,p} \times (G^{[p]}(\mathbb{R}^d)))$, which, at the first level is the classical support theorem of Stroock–Varadhan [25].

If $y_0$ and the $V_i$'s are the image by a continuous function of $\mathbf{B}$, then the support is still trivially characterized. One could then ask for more specific conditions on $y_0$ and the $V_i$'s, in the spirit of [17], and obtain a detailed support theorem. Thanks to the universal limit, one would obtain stronger results than in [17] (stronger topology and without the assumption of deterministic vector fields) but we shall not pursue this here.

L. COUTIN
LABORATOIRE DE STATISTIQUE ET PROBABILITÉS
UNIVERSITÉ PAUL SABATIER
118 ROUTE DE NARBONNE
31062 TOULOUSE CEDEX
FRANCE
E-MAIL: coutin@cict.fr

P. FRIZ
CENTER FOR MATHEMATICAL SCIENCES
UNIVERSITY OF CAMBRIDGE
WILBERFORCE ROAD
CAMBRIDGE CB3 0WB
UK
E-MAIL: P.Friz@statslab.cam.ac.uk

N. VICTOIR
UNIVERSITY OF OXFORD
MAGDALEN COLLEGE
OXFORD OX1 4AU
UK
E-MAIL: Victoir@maths.ox.ac.uk